\documentclass[conference, 9pt]{IEEEtran}
\IEEEoverridecommandlockouts

\usepackage{cite}
\usepackage{amsmath,amssymb,amsfonts}
\usepackage{algorithm}
\usepackage{algorithmic}
\usepackage{graphicx}
\usepackage{textcomp}
\usepackage{xcolor}
\usepackage{enumitem}
\def\BibTeX{{\rm B\kern-.05em{\sc i\kern-.025em b}\kern-.08em
    T\kern-.1667em\lower.7ex\hbox{E}\kern-.125emX}}

\usepackage{amssymb}
\usepackage{amsthm}

\usepackage{mathtools}

\usepackage{hyperref}
\hypersetup{
  colorlinks   = true, 
  urlcolor     = blue, 
  linkcolor    = blue, 
  citecolor    = blue   
}
\usepackage{bm}
\usepackage{tabularx}
\usepackage{pifont}

\newtheorem{theorem}{Theorem}[section]

\newtheorem{lemma}[theorem]{Lemma}

\newtheorem{assumption}[theorem]{Assumption}

\newtheorem{remark}[theorem]{Remark}

\def\Prm{\mathbf{X}}
\def\prm{\mathbf{x}}
\def\avgprm{\bar{\mathbf{X}}}
\def\dprm{\bm{\lambda}}

\def\wta{\widetilde{\bf A}}

\def\wtq{\widetilde{\bf Q}}
\def\wtk{\widetilde{\bf K}}

\def\algoname{TiCoPD}

\def\m{M}

\def\serr{{\bf e}_s}

\def\gerr{{\bf e}_g}

\def\wxv{{\bf W}_{xv}}

\def\a{{\tt a}}
\def\b{{\tt b}}
\def\c{{\tt c}}
\def\d{{\tt d}}
\def\e{{\tt e}}

\newcommand{\dotp}[2]{\left\langle{#1}\ \middle|\ {#2}\right\rangle}
\newcommand{\expec}[1]{\mathbb{E}\left[ {#1} \right]}

\newcommand*{\htwai}[1]{\textbf{\textcolor{red}{To: #1}}}

\begin{document}

\title{A Two-timescale Primal-dual Algorithm for Decentralized Optimization with Compression
}

\author{\IEEEauthorblockN{Haoming Liu, Chung-Yiu Yau, Hoi-To Wai}
\IEEEauthorblockA{\emph{Dept. of Systems Engineering \& Engineering Management}, \textit{The Chinese University of Hong Kong}\\ Emails:~\url{haomingliu@cuhk.edu.hk}, \url{cyyau@se.cuhk.edu.hk}, \url{htwai@se.cuhk.edu.hk}.} 
\thanks{This research is supported in part by project \#MMT-p5-23 of the Shun Hing Institute of Advanced Engineering, The Chinese University of Hong Kong.}}


\maketitle

\begin{abstract}
This paper proposes a two-timescale compressed primal-dual (TiCoPD) algorithm for decentralized optimization with improved communication efficiency over prior works on primal-dual decentralized optimization.  
The algorithm is built upon the primal-dual optimization framework and utilizes a majorization-minimization procedure. The latter naturally suggests the agents to share a compressed difference term during the iteration. 
Furthermore, the TiCoPD algorithm incorporates a fast timescale mirror sequence for agent consensus on nonlinearly compressed terms, together with a slow timescale primal-dual recursion for optimizing the objective function.
We show that the TiCoPD algorithm converges with a constant step size. It also finds an $\mathrm{\mathcal{O}( {\it 1 /T} )}$ stationary solution after $\mathrm{\it T}$ iterations. Numerical experiments on decentralized training of a neural network validate the efficacy of TiCoPD algorithm.
\end{abstract}


\section{Introduction}
Let ${\cal G} = (V,E)$ be an undirected and connected graph of $n$ agents, with the node set given by $V = [n] := \{1,...,n\}$ and the edge set $E \subseteq V \times V$. Consider the following distributed optimization problem:
\begin{equation} 
    \min_{ \Prm \in \mathbb{R}^{nd}} \quad \frac{1}{n} \sum_{i=1}^n f_i( \Prm_i ) \quad \text{s.t.} \quad \Prm_i = \Prm_j, ~\forall (i,j) \in E, \label{eq:main_problem}
\end{equation}
where $\Prm = [ \Prm_1^\top ~ \Prm_2^\top ~\cdots ~ \Prm_n^\top]^\top$. For each $i \in V$, $\Prm_i \in \mathbb{R}^d$ stands for the local decision variable of the $i$th agent, and the local objective function $f_i: \mathbb{R}^d \rightarrow \mathbb{R}$ is continuously differentiable (possibly non-convex).
Distributed solution methods for \eqref{eq:main_problem} has found applications in machine learning \cite{mateos2010distributed, lian2017can}, signal processing \cite{dimakis2010gossip, xiao2007distributed, kar2013consensus}, etc. 


In scenarios where the local objective function (and the associated training data) has to be stored or processed in a distributed network, and/or there is no central server, a decentralized algorithm is preferred for tackling \eqref{eq:main_problem} where agents do not share their local objective functions with other parties. The pioneering work by Nedi\'{c} and Ozdaglar \cite{nedic2009distributed} proposed the decentralized gradient descent (DGD) method. The latter is shown to converge in the non-convex setting in a follow-up publication \cite{zeng2018nonconvex}. A number of subsequent works have been proposed to improve DGD, e.g., EXTRA in \cite{shi2015extra}, gradient tracking in \cite{qu2017harnessing}, etc. The (proximal) primal-dual algorithm \cite{chang2020distributed, hajinezhad2019perturbed, yi2021linear} yields a general framework for decentralized optimization with good convergence property where it encompasses the previous algorithms as special cases.


For high dimensional instances of \eqref{eq:main_problem} with $d \gg 1$, the communication overhead incurred with decentralized optimization can be a bottleneck impeding the convergence speed of these algorithms. Reducing the bandwidth usage through lossy communication compression has thus become an important issue in recent works. Note that directly applying compression scheme such as quantization in distributed optimization may result in a non-converging algorithm \cite{rabbat2005quantized}. Thus, many works have proposed to combine the decentralized algorithms with an error feedback subroutine to achieve exact convergence. For example, \cite{reisizadeh2019exact, magnusson2020maintaining, liu2021linear, zhao2022beer} focused on the deterministic gradient setting and extended algorithms such as DGD and gradient tracking with communication compression, \cite{koloskova2019decentralized, koloskova2019decentralizedb, yau2022docom, xie2024communication} considered similar extensions but have considered using stochastic gradients. Alternatively, \cite{stich2018local, yu2019parallel, basu2020qsparse} considered reducing the communication frequency by local updates, \cite{singh2019sparq} considered event triggered communication, \cite{michelusi2022finite, nassif2024differential} considered algorithms that employ an adaptive finite-bit quantizer.

Meanwhile, existing results entail various limitations in guaranteeing convergence for the general non-convex setting --- the algorithms in \cite{reisizadeh2019exact, koloskova2019decentralizedb} require a bounded heterogeneity assumption such that $\|\nabla f_i(\Prm)-\nabla f_j(\Prm)\|$ is upper bounded for any $\Prm$ and may require a diminishing stepsize for convergence, the algorithms in \cite{zhao2022beer, yau2022docom} requires applying two separate communication compression operations that may increase storage complexity, and other algorithms \cite{koloskova2019decentralized, magnusson2020maintaining, liu2021linear} have not been analyzed for the non-convex optimization setting. 

This paper aims to address the above shortcomings with communication efficient decentralized optimization. Our key idea is to develop the algorithm using the primal-dual framework \cite{hajinezhad2019perturbed, chang2020distributed} and to incorporate compressed communication updates under the classical majorization-minimization framework \cite{sun2016majorization}. 
Our contributions are:
	\begin{itemize}[leftmargin=*]
		\item We propose the {\bf T}wo-t{\bf i}mescale {\bf Co}mpressed {\bf P}rimal-{\bf D}ual ({\algoname}) algorithm as a nonlinearly compressed primal-dual algorithm for decentralized optimization. {\algoname} follows a two-timescale update rule which separates the \emph{communication} and \emph{optimization} steps, handled using a large stepsize and a small stepsize, respectively. 
		\item To incorporate nonlinear compression into the primal-dual algorithm, we develop a majorization-minimization (MM) procedure which suggests agents to transmit the compressed difference terms -- a scheme that coincides with the popular error feedback mechanism. This offers a new perspective for extending the error feedback mechanism popularized by \cite{koloskova2019decentralized} and establishing a connection to the nonlinear gossiping algorithm in \cite{mathkar2016nonlinear}. We believe that this observation will be of independent interest. 
		\item For optimization problems with continuously differentiable (possibly non-convex) objective functions, we show that {\algoname} converges at a rate of ${\cal O}(1/T)$ towards a stationary solution of \eqref{eq:main_problem}. Our result is proven with a constant stepsize and without relying on additional assumptions such as bounded heterogeneity.
	\end{itemize}
Finally, we present numerical experiments to demonstrate the efficacy of {\algoname} algorithm against state-of-the-art algorithms in tackling a toy example of training neural networks over network for \eqref{eq:main_problem}.


\section{Problem Statement} \label{sec:Problem Statement}
This section introduces the basic ideas for the proposed {\algoname} algorithm to tackling \eqref{eq:main_problem}. We first introduce a few extra notations to facilitate the development: define
    \begin{align} 
        & \textstyle f ( \avgprm ) := \frac{1}{n} \sum_{i=1}^n f_i( \avgprm ),\quad \nabla f( \avgprm) := \frac{1}{n} \sum_{i=1}^n \nabla f_i( \avgprm ),
    \end{align}
as the global objective function and gradient evaluated on a common decision variable $\avgprm \in \mathbb{R}^d$, respectively. Moreover, the consensus constraint $\Prm_i = \Prm_j, ~\forall (i,j) \in E$ can be replaced by the equality $\wta \Prm = {\bf 0}$, where $\wta = {\bf A}\otimes {\bf I}_d$ and $\wta \in \{0, 1, -1\}^{|E| \times n}$ denotes the incidence matrix of the graph $G$.

Using the above notations and introducing the Lagrange multiplier variable $\dprm = ( \dprm_i )_{i \in E} \in \mathbb{R}^{ |E| d }$ for the equality constraint $\wta \Prm = {\bf 0}$, we consider the following augmented Lagrangian function:
	\begin{equation} \label{eq:aug_lag}
		\mathcal{L}( \Prm, \dprm) := \frac{1}{n} \sum_{i=1}^n f_i( \Prm_i  ) + \dprm^\top \wta \Prm + \frac{\theta}{2} \|\wta \Prm\|^2 ,
	\end{equation} 
where $\theta > 0$ is a regularization parameter. It can be shown that any stationary point to ${\cal L}(\cdot)$ satisfying $\nabla_{ \Prm } \mathcal{L}( \Prm, \dprm) = {\bm 0}$, $\nabla_{ \dprm } \mathcal{L}( \Prm, \dprm) = {\bm 0}$ is a KKT point of \eqref{eq:main_problem}.
Importantly, applying the standard gradient descent-ascent (GDA) algorithm on ${\cal L}(\cdot)$ yields the decentralized algorithm proposed in \cite{yi2021linear} which finds a stationary solution to \eqref{eq:main_problem} for a class of non-convex problems.

However, a key drawback for \cite{yi2021linear} and similar algorithms proposed using the primal-dual optimization framework (e.g., \cite{hajinezhad2019perturbed}) lies on their high bandwidth usage for high-dimensional problems when $d \gg 1$. In particular, at each iteration, such algorithms require the agents to share their current local decision variable with neighboring agents, which demands transmitting an $\mathbb{R}^d$ vector. Such step introduces a considerable communication overhead for their implementation.

\section{Proposed TiCoPD Algorithm} \label{sec:2TS Compressed Stochastic PDA}
Our idea is to develop a {\bf T}wo-t{\bf i}mescale {\bf Co}mpressed {\bf P}rimal-{\bf D}ual ({\algoname}) algorithm that supports general \emph{nonlinear compression} in primal-dual decentralized updates. Unlike prior works such as CHOCO-SGD \cite{koloskova2019decentralized} that takes the compressed update as a separate sub-routine, we treat \emph{compressed communication} and \emph{optimization} as lower and upper level updates, respectively, to be updated at different speed. The algorithm depends on two ingredients: (i) a majorization-minimization step that introduces a \emph{surrogate} variable to separate the communication step from the optimization step, (ii) a \emph{two-timescale} update that incorporates the nonlinearly compressed update of the surrogate variable. The above ideas will be introduced in the sequel. 

\vspace{.1cm}
\noindent {\bf Majorization Minimization.} We begin by inspecting the augmented Lagrangian function ${\cal L}(\cdot)$ again. Here, computing the gradient for the last term in \eqref{eq:aug_lag} at a fixed primal variable $\Prm^t$ leads to 
\begin{equation} \textstyle 
\nabla_{ \Prm_i } || \wta \Prm^t ||^2 = 2 \sum_{j \in {\cal N}_i } ( \Prm_j^t - \Prm_i^t ),
\end{equation}
where ${\cal N}_i$ is the neighbor set of agent $i$.
Consequently, the primal update necessitates the communication of the neighbors' decision variables $\Prm_j$, $j \in {\cal N}_i$ and leads to a communication bottleneck. 

Our idea is to sidestep this term through a majorization-minimization (MM) procedure with a surrogate variable. 
Let there be a sequence of surrogate variables $\{ \hat\Prm^t \}_{t \geq 0}$ such that (i) $\hat{\Prm}^t \approx \Prm^t$, and (ii) it is possible for agent $i$ to acquire the neighbors' surrogate variables $( \hat{\Prm}_j^t )_{j \in {\cal N}_i^t}$ with \emph{compressed communication}. We will illustrate how to construct such sequence later. With $\m := \| \wta^\top \wta \|_2$, the following \emph{majorization} holds for any $\Prm$,
\begin{equation}
	\begin{aligned}
	& \| \wta \Prm \|^2 \le  \| \wta\hat{\Prm}^t\|^2 +2 ( \Prm -\hat{\Prm}^t)^\top \wta^\top \wta\hat{\Prm}^t  +\m \| \Prm -\hat{\Prm}^t\|^2 .
	\end{aligned}
\end{equation}
Unlike the original term $\| \wta \Prm \|^2$, evaluating the gradient w.r.t.~$\Prm_i$ on the upper bound only requires aggregating the surrogate variables $\hat{\Prm}_j^t, j \in {\cal N}_i$, which can be obtained from compressed communication.

We further upper bounding the first term in \eqref{eq:aug_lag} using standard truncated Taylor approximation. Now, the $\Prm$-update can be computed using the following \emph{minimization} step: 
\begin{align} 
	&\Prm^{t+1} =  \arg\min_{\Prm \in \mathbb{R}^{nd} } \nabla {\bf f}( \Prm^t)^\top (\Prm - \Prm^t) + \Prm^\top \wta^\top \dprm^t +\frac{\theta}{2} \| \wta\hat{\Prm}^t\|^2 \notag \\
	& \quad \qquad \qquad +\theta \Prm^\top \wta^\top \wta\hat{\Prm}^t +\frac{\theta \m}{2}\| \Prm-\hat{\Prm}^t\|^2 +\frac{1}{2\widetilde{\alpha}}\| \Prm -\Prm^t\|^2 \notag \\
    & = \beta \Prm^t + (1-\beta) \hat{\Prm}^t - \alpha \big(\nabla {\bf f}(\Prm^t)+ \wta^\top {\dprm}^{t} +\theta \wta^\top \wta\hat{\Prm}^t \big), \label{eq:x-update1}
\end{align}
where $\nabla {\bf f}( \Prm^t) = [ \nabla f_1( \Prm_1^t)^\top \cdots \nabla f_n( \Prm_n^t )^\top ]^\top$, $\alpha = \frac{1}{\frac{1}{\widetilde{\alpha}}+\theta \m}$, and $\beta = \frac{\alpha}{\widetilde{\alpha}}$.
For the $\dprm$-subproblem, similarly we replace $\Prm^t$ with the surrogate variable $\hat{\Prm}^t$ to obtain:
\begin{equation} \textstyle \label{eq:lambda-update1} 
\begin{aligned}	
 \dprm^{t+1} & = \arg\min_{\dprm \in \mathbb{R}^{|E|d}} ~ \left\{ -\dprm^\top \wta\hat{\Prm}^t +\frac{1}{2\eta} \| \dprm-\dprm^t \|^2 \right\} \\
 & = \dprm^t + \eta \wta \hat{\Prm}^t .
\end{aligned}
\end{equation}
Substituting the variable $\widetilde{\dprm}^t = \wta^\top \dprm^t \in \mathbb{R}^{nd}$ yields a decentralized primal-dual algorithm that only requires aggregating $\hat{\Prm}_i^t$ at each step.
Our remaining task is to study how to effectively construct the sequence $\{ \hat{\Prm}^t \}_{t \geq 0}$ in a compressed communication friendly fashion.

\vspace{.1cm}
\noindent {\bf Two-timescale Updates.}
Recall that a key requirement for our construction of $\{ \hat{\Prm}^t \}_{t \geq 0}$ is that the surrogate variable should \emph{track} the original decision variable $\hat{\Prm}^t \approx \Prm^t$. 

To this end, our construction consists of applying a (randomized) {compression operator} $\hat{Q}: \mathbb{R}^d \times \Omega_i \to \mathbb{R}^d$ satisfying: there exists $\delta \in (0,1]$ such that
\begin{equation} \label{eq:contractive}
    \expec{ \| \hat{Q} (\prm; \xi_{i})-\prm\|^2 }  \le (1-\delta)^2 \expec{\|\mathbf{\prm}\|^2 },~\forall~ \prm \in \mathbb{R}^d, i\in[n],
\end{equation}
where $\xi_i \in \Omega_i$ is a random variable.
For example, the above property can be satisfied with the randomized quantization operator:
\begin{equation} \label{eq:qsgd}
	\mathrm{qsgd}_s( \prm ; \xi_{i})=\frac{\mathrm{sign}(\prm)\cdot\|\prm\|}{s\tau}\cdot\left\lfloor s\frac{|\prm|}{\|\prm\|}+\xi_{i} \right\rfloor\mathrm{~,}
\end{equation}
where $s > 0$ is the number of precision levels, $\tau = 1 + \min\{ d/s^2, \sqrt{d}/s \}$ and $\xi_{i} \sim {\cal U}[0,1]^d$ is an additive noise \cite{alistarh2017qsgd}. The operator can be implemented with an \emph{encoder} (denoted ${\rm ENC}(\cdot)$) which turns the $d$-dimensional input vector into a $d \log_2 s$ bits string to represent the quantized levels, a $d$ bits string to represent $\mathrm{sign}(x)$ and a floating-point value for $\|\prm\|$. The latter can be transmitted on bandwidth limited channels. On the receiver's side, a \emph{decoder}  (denoted ${\rm DEC}(\cdot)$) may convert the received symbols into the quantized real vector in \eqref{eq:qsgd}. Particularly, the randomized quantization operator satisfies \eqref{eq:contractive} with $\delta = \frac{1}{2 \tau}$; see \cite{alistarh2017qsgd, koloskova2019decentralized} for discussions on contractive compressors satisfying \eqref{eq:contractive}. Finally, we denote ${Q}:\mathbb{R}^{nd} \times \Omega_1\times\ldots\times\Omega_n \to \mathbb{R}^{nd}$ as the compression operator for the column stacked variables such that $[{Q}(x;\xi)]_i = \hat{Q}(x_i;\xi_{i})$, where $\xi = (\xi_1,\ldots,\xi_n)$.

The contraction property in \eqref{eq:contractive} suggests a compressed communication friendly procedure to achieve $\hat{\Prm}^t \approx \Prm^t$. In particular, we notice that the mean field $\mathbb{E}_{\xi} [ Q(\Prm^{t} - {\Prm} ;\xi )]$ has a unique fixed point at $\Prm = \Prm^t$. In particular, we observe that at the $t$th primal-dual iteration, if we let $k$ denotes the contraction iteration index and $\gamma \in (0,1]$ be a stepsize parameter, the recursion
\begin{equation} \label{eq:xhat_update}
	\hat{\Prm}^{t,k+1} = \hat{\Prm}^{t,k} + \gamma \, Q(\Prm^{t} - \hat{\Prm}^{t,k} ;\xi^{t,k+1}),~\forall~k \geq 0.
\end{equation}
finds $\hat{\Prm}^{t,k} \overset{k \to \infty}{\to} \Prm^t$.
The insight behind the above update is that the compression error can be gradually reduced by compressing only the error itself $\Prm^t - \hat{\Prm}^{t,k}$.
For implementation on a decentralized system, for each $(t,k)$, the agents only need to encode and transmit the differences $\Prm^t - \hat{\Prm}^{t,k}$. In particular, it holds $\mathbb{E} [ \| \hat{\Prm}^{t,k} - \Prm^t \|^2 ] \leq (1-\gamma \delta)^{k} \| \hat{\Prm}^{t,0} - \Prm^t \|^2 $.
For the rest of this paper, we shall take $\gamma = 1$ for simplicity. However, we remark that in cases when \eqref{eq:contractive} is not satisfied, the stepsize parameter $\gamma \in (0,1]$ can be used to control the process for the construction of $\hat{\Prm}^t \approx \Prm^t$.

The convergence of the procedure \eqref{eq:xhat_update} relies on holding $\Prm^t$ fixed as $k$ increases. While repeating the recursion \eqref{eq:xhat_update} until $k \to \infty$ guarantees $\hat{\Prm}^{t,k} = \Prm^t$, this may incur a prohibitive communication cost for each of the primal-dual update \eqref{eq:x-update1}, \eqref{eq:lambda-update1} indexed by $t$. Fortunately, we observe that each recursion in \eqref{eq:xhat_update} can actually reduces the error by a factor of $(1- \delta)$, i.e., it exhibits a geometric convergence rate.
Suppose the primal-dual update stepsizes are sufficiently small compared to $\gamma$, it is possible to satisfy the \emph{tracking condition} $\hat{\Prm}^t \approx \Prm^t$ with only one step of the update in \eqref{eq:xhat_update} per primal-dual iteration.  

Overall, the above discussion suggests a \emph{two-timescale} simultaneous update of \eqref{eq:x-update1}, \eqref{eq:lambda-update1}, \eqref{eq:xhat_update}. We call this algorithm the {\algoname} algorithm which is summarized in Algorithm~\ref{alg:TiCoPD}.
Note that the algorithm is fully decentralized. For each iteration, the agents always communicate with one round of \emph{compressed} message exchanges.

\begin{algorithm}[t] 
    \caption{\algoname} 
    \begin{algorithmic}[1]\label{alg:TiCoPD}
    \STATE {\bf Input:} Parameter $\alpha,\theta,\beta,\eta$, 
    initialize $\Prm^0,\hat{\Prm}^0,\widetilde{\lambda}^0$.
    \FOR{$t = 1,\ldots,T$}
        \STATE \emph{Surrogate variable update:} for any $i \in V$,
            \begin{equation} \notag
                \hat{\Prm}^{t}_i = \hat{\Prm}^{t-1}_i + \hat{Q}(\Prm^{t}_i - \hat{\Prm}^{t-1}_i;\xi_{i}^{t}),
            \end{equation}
        and transmit the compressed message ${\rm ENC}( \Prm^{t}_i - \hat{\Prm}^{t-1}_i ; \xi_i^t )$ to agent $j \in {\cal N}_i$. 
        \STATE \emph{Aggregate received messages:} for any $i \in V$,
            \[ \textstyle 
                \hat{\Prm}_{i, -i}^{t} = \hat{\Prm}_{i,-i}^{t-1} +  \sum_{j \in {\cal N}_i} {\rm DEC} ( {\rm ENC} (\Prm^{t}_j - \hat{\Prm}_{j}^{t-1};\xi_{j}^{t}) ), \vspace{-.0cm}
            \]
            where $\hat{Q}(\cdot) = {\rm DEC} \circ {\rm ENC} (\cdot)$ and ${\rm DEC}(\cdot)$ is a decoder.
            \STATE \emph{Primal-dual update:} for any $i \in V$,
            \begin{align}
                \Prm^{t+1}_i & = \beta \Prm^t_i + (1-\beta) \hat{\Prm}^t_i \notag \\
                & \quad - \alpha \big[ \nabla f_i(\Prm_i^t) + \widetilde{\dprm}_i^{t}+ \theta ( |{\cal N}_i| \hat{\Prm}^t_i - \hat{\Prm}^t_{i,-i} ) )\big] , \notag \\
                \widetilde{\dprm}^{t+1}_i & \textstyle = \widetilde{\dprm}^t_i + \eta ( |{\cal N}_i| \hat{\Prm}^t_i - \hat{\Prm}^t_{i,-i} ). \notag \vspace{-.2cm}
            \end{align}
    \ENDFOR
    \end{algorithmic}
    \end{algorithm}


\begin{remark}
The recently proposed CP-SGD algorithm \cite{xie2024communication} uses a similar idea of contractive compressor on primal-dual algorithm as {\algoname}. There are several key differences: (i) CP-SGD uses an extra auxiliary variable $\Prm^c$ and entails an additional stepsize parameter $\alpha_x$, (ii) in practice we observe that CP-SGD yields slower convergence in the consensus error for certain problems. We highlight that {\algoname} was developed directly from the MM procedure and two-timescale updates. Our algorithm also belongs to a general framework that can incorporate scenarios with less restrictions on the communication architecture, e.g., noisy communication, time varying graphs, etc. Such features are missing in the CP-SGD algorithm \cite{xie2024communication}.
\end{remark}

\section{Convergence analysis}  \label{sec:Convergence analysis}
This section establishes the convergence of the {\algoname} algorithm towards a stationary point of \eqref{eq:main_problem} at a sublinear rate. To facilitate our discussions, we define the consensus error operator $\wtk := ( {\bf I}_n - {\bf 1}{\bf 1}^\top/n ) \otimes {\bf I}_d$ and $\wtq := ( \wta^\top \wta )^\dagger$, where $(\cdot)^\dagger$ denotes the Moore-Penrose inverse. From the definitions, we observe that $\wta^\top \wta \wtk = \wta^\top \wta = \wtk \wta^\top \wta$ and $\wtq \wta^\top \wta = \wta^\top \wta \wtq = \wtk$. We first state the assumptions about the objective function, graph, and the compression operator. 
    
\begin{assumption} \label{assm:lip}
For any $i \in [n]$, the function $f_i$ is $L$-smooth, i.e.,
\begin{equation} \label{eq:f_lip} \textstyle
\|\nabla f_i({\bf x}) - \nabla f_i({\bf y}) \| \leq L \|{\bf x} - {\bf y} \| ~\forall~ {\bf x},{\bf y} \in \mathbb{R}^{d}.
\end{equation}
\end{assumption}
\noindent Note that \eqref{eq:f_lip} implies the global objective function is $L$-smooth. It is also a standard assumption in the optimization literature. 
	
	
\begin{assumption} \label{assm:rand-graph}
There exists $\tilde{\rho}_1 \ge \tilde{\rho}_2> 0$ such that 
\begin{equation} \textstyle
\begin{aligned} 
& \tilde{\rho}_2 \wtk \preceq \wta^\top \wta \preceq \tilde{\rho}_1 \wtk. \label{eq:q_ineq}
\end{aligned}
\end{equation}
\end{assumption}
\noindent Recall that as $\wta^\top \wta$ is the graph Laplacian matrix, the assumption above is satisfied when $G$ is a connected graph. Furthermore, we have $\tilde{\rho}_1 = \| \wta^\top \wta \|_2$. It also holds $\tilde{\rho}_1^{-1} \wtk \preceq \wtq \preceq \tilde{\rho}_2^{-1} \wtk$. 
	
\begin{assumption} \label{assm:compress}
For any fixed $\Prm \in \mathbb{R}^d$, there exist $0 < \delta\le 1$ such that the compression operator $\hat Q(\cdot; \xi_q)$ satisfies
\begin{equation} \textstyle
\expec{ \| \hat Q(\Prm; \xi_q)-\Prm\|^2 }  \le (1-\delta)^2  \| \Prm \|^2. \label{eq:Q-contraction}
\end{equation}
\end{assumption}
\noindent As discussed before, the randomized quantizer \eqref{eq:qsgd} satisfies the above with $\delta = s/\sqrt{d}$ for sufficiently large $d$. Other common compression operators like  random-k sparsification, top-k sparsification also satisfy the assumption, as shown in \cite{koloskova2019decentralized}.

To describe our main results, we need to define a few extra notations. The network-average decision variable is denoted as
\begin{equation} \label{eq:avgprm} 
\avgprm^t = \frac{1}{n} \sum_{i=1}^n \Prm_i^t = \frac{1}{n} ( {\bf 1}^\top \otimes {\bf I}_d ) \Prm^t, \notag
\end{equation}
and the consensus error is given by
\begin{equation} 
\sum_{i=1}^n \| \Prm_i^t - \avgprm^t \|^2 = \| \Prm^t - {\bf 1}_n \otimes \avgprm^t \|^2 = \| \wtk \Prm^t \|^2 = \| \Prm^t \|_{\wtk}^2. \notag
\end{equation}
We also set $\mathbf{v}^t = \alpha ( \widetilde{\dprm}^{t} + \nabla {\bf f}(({\bf 1} \otimes {\bf I}_d ) \avgprm^{t}) )$ 
as an auxiliary variable to measures the violation of tracking the average deterministic gradient. 
The convergence result follows:

\begin{theorem} \label{theo:main_theorem}
Under Assumptions \ref{assm:lip}--\ref{assm:compress}, suppose the step sizes satisfy $\eta >0, \theta \ge \theta_{lb}, \alpha\le\alpha_{ub}$, where
\begin{equation}
\begin{aligned}
\theta_{lb} = & \max \Big\{\frac{4L^2}{n \tilde{\rho}_2 \a}, \frac{2}{ \tilde{\rho}_2} (1+2L+8\eta\tilde{\rho}_1^2\tilde{\rho}_2^{-1} \\ 
& \qquad + \delta_1(\frac{3}{2}+3L^2+\eta\tilde{\rho}_1+\frac{\tilde{\rho}_1^2}{2})) \Big\}, \\
\alpha_{ub} =& \frac{\max\{ 16\eta, \delta\} }{320\theta^2} \cdot \min\Big\{\frac{1}{\m^2},\frac{n\a}{\m^2},\frac{1}{\tilde{\rho}_1^{2}}\Big\},
\end{aligned}
\end{equation}
where $\delta_2 = \max\{\frac{16 \eta}{\delta},1 \}$, $\delta_1 = 12 \max\{ 2, 2\tilde{\rho}_2^{-1} \eta^{-1}, \delta_2\tilde{\delta}\}$, $\tilde{\delta}=\max\{\frac{(1-\delta)^2(1-\frac{\delta}{2})^2}{(1-\frac{\delta}{2})^2-(1-\delta)^2},1\}$.
Then, for any $T \geq 1$, it holds  
\begin{equation}
\begin{aligned}
& \frac{1}{T} \sum_{t=0}^{T-1}\expec{\left\| \nabla f (\avgprm^{t}) \right\|^2} \leq \frac{ \expec{F_0} - \expec{F_{T}}}{\alpha T/16},\\
&\frac{1}{T} \sum_{t=0}^{T-1}\expec{ \| \Prm^t \|_{\wtk}^2} \le  \frac{ \expec{F_0} - \expec{F_{T}}}{\alpha\theta\tilde{\rho}_2 \a T/4} .
\end{aligned}
\end{equation}
where the expectation is taken w.r.t.~randomness in the compression operator, $\a>0$ is a free quantity, and
\begin{equation}
\begin{split}
F_{t} & = f (\avgprm^{t}) - f^\star +
\frac{\a}{\eta \alpha} \|\mathbf{v}^{t}\|_{{\widetilde{\mathbf{Q}}+
\alpha \left( (\frac{\delta_1}{2}(\theta+\eta)-\theta)-\delta_2 \theta \right) \widetilde{\mathbf{K}}}}^{2} \\
& +
\a \|\Prm^{t}\|_{{\widetilde{\mathbf{K}}}}^{2} + \delta_1 \a \left\langle\Prm^{t}\mid\mathbf{v}^{t}\right\rangle_{{\widetilde{\mathbf{K}}}} +
\delta_2 \a \| \hat{\Prm}^t-\Prm^t \|^2,
\end{split}
\end{equation}
can be shown to be non-negative.
\end{theorem}
\noindent Our analysis involves constructing a Lyapunov function to track the joint convergence of the coupled variables. However, the proof is rather technical and are relegated to an online appendix; see \url{https://www1.se.cuhk.edu.hk/~htwai/pdf/icassp25-ticopd.pdf}.

The above theorem shows that upon fixing the dual step size at $\eta=\delta$, there exists a sufficiently small constant primal step size $\alpha \leq \alpha_{ub}$\footnote{Since it holds $\alpha=\frac{1}{\frac{1}{\widetilde{\alpha}}+\theta L}$, once $\theta$ is fixed, the upper bound of $\alpha$ can be achieved by selecting a sufficiently small $\widetilde{\alpha}$.} (independent of $T$) such that for sufficiently large $T$,
\begin{align}
\frac{1}{T} \sum_{t=0}^{T-1}\expec{\left\| \nabla f (\avgprm^{t}) \right\|^2} = {\cal O} \big( 1 / T \big), 
\end{align}
where it implies that there exists $\hat{t} \in \{0, \ldots,T-1\}$ such that $\expec{ \| \nabla f (\avgprm^{\hat{t}}) \|^2} = {\cal O}(1/T)$.
Notice that this is the same rate (w.r.t.~$1/T$) as a centralized gradient algorithm on the smooth optimization problem \eqref{eq:main_problem}. Moreover, compared to the results for decentralized algorithms such as DGD \cite{zeng2018nonconvex} and CHOCO-SGD \cite{koloskova2019decentralizedb}, our results do not require using a diminishing step size nor bounded gradient heterogeneity. 

We also comment on the impact of compression on the convergence rate and step size selection. Recall that $\delta \in (0,1]$ of Assumption~\ref{assm:compress} is affected by quality of the compressor, where $\delta \approx 0$ with an aggressive compression scheme, e.g., the number of quantization levels, $s$, is small. In this case, we have $\tilde{\delta} \asymp \delta^{-1}$, $\theta_{lb} \asymp \delta_1 \asymp \delta^{-1}$, and thus $\alpha_{ub} = {\cal O} (\delta^3)$. Consequently, we observe that the convergence rate will be $\frac{1}{T} \sum_{t=0}^{T-1}\expec{\left\| \nabla f (\avgprm^{t}) \right\|^2} = {\cal O} ( \delta^3 T^{-1} )$ and the upper bound evaluates to ${\cal O} ( d^{1.5} s^{-3} T^{-1} )$ for the case of randomized quantization. As such, we conclude that the convergence of {\algoname} is sensitive to the quality of compressor. A future direction is to improve such dependence on $\delta$.



\section{Numerical Experiments} \label{sec:Numerical Experiments}
We compare the performance of the proposed algorithm on a typical machine learning model training task. Here, $f_i( \Prm_i )$ is taken as the cross-entropy classification loss where $\Prm_i \in \mathbb{R}^{79510}$ denotes the weights of a 2 layer feed-forward neural network with 100 neurons and sigmoid activation. We consider training the neural network on the MNIST dataset. We benchmark the performance among algorithms that support compressed communication including CHOCO-SGD \cite{koloskova2019decentralized} and CP-SGD \cite{xie2024communication} (both implemented with exact gradients), while also comparing against DEF-ATC \cite{nassif2024differential}, BEER \cite{zhao2022beer}, the classical DGD method, and a heuristic variant of DGD that only aggregates quantized parameters. 
The stepsizes for each of the above algorithms are fine tuned so that they achieve the best performance after $10^5$ iterations.
The experiments are performed on a 40-core Intel Xeon server with 64GB memory using the PyTorch package.

We simulate the algorithms on a ring network of 10 nodes. The MNIST dataset (of $M = 60,000$ images, each with $p=784$ pixels) is distributed to each node according to the class label $\{0,...,9\}$, i.e., each node only access one class of images. We consider unshuffled dataset to maximize data heterogeneity across agents. 
Note that the convergence of primal-dual algorithms ({\algoname}, CP-SGD) are unaffected by the data heterogeneity issue. 

\begin{figure}
    \centering
    \includegraphics[width=0.45\textwidth]{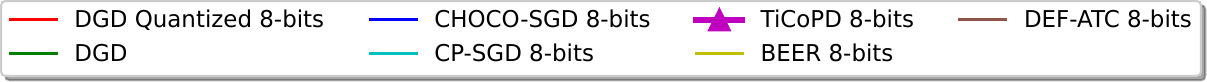} \\
    \vspace{0.25cm}
    \includegraphics[width=0.235\textwidth]{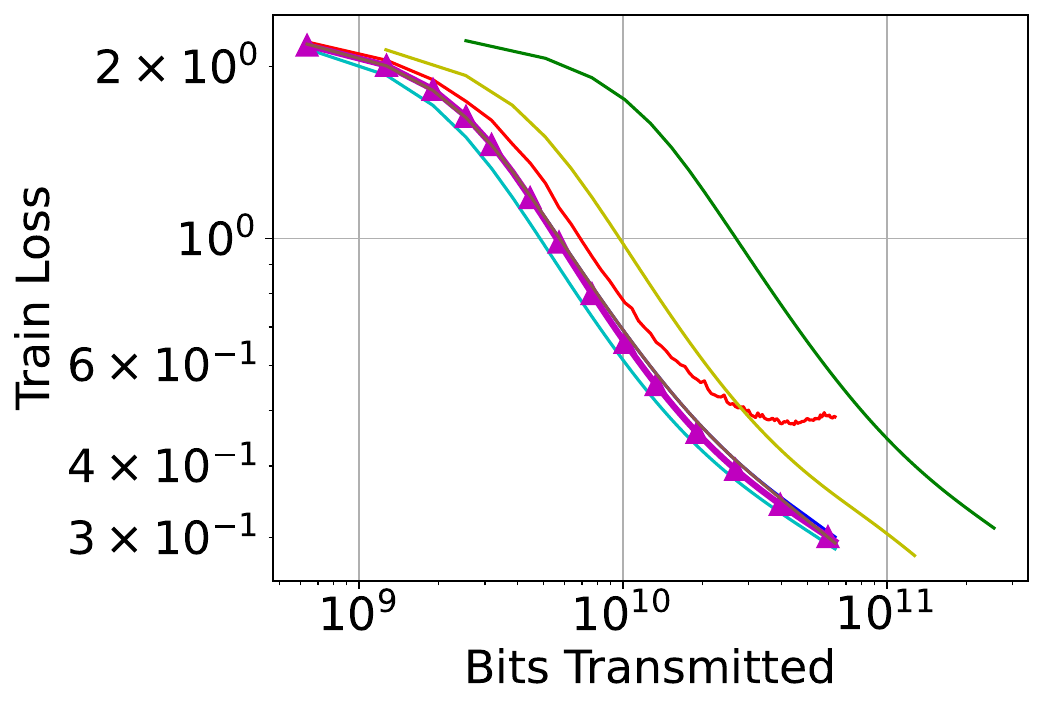}
    \includegraphics[width=0.235\textwidth]{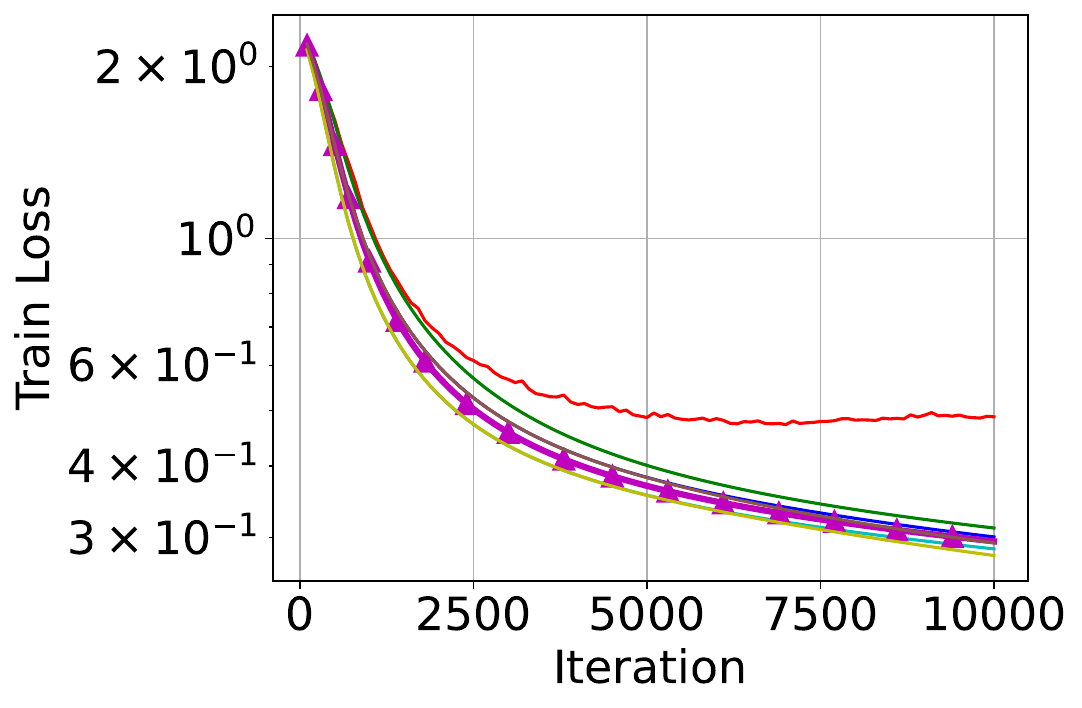} \\
    \includegraphics[width=0.235\textwidth]{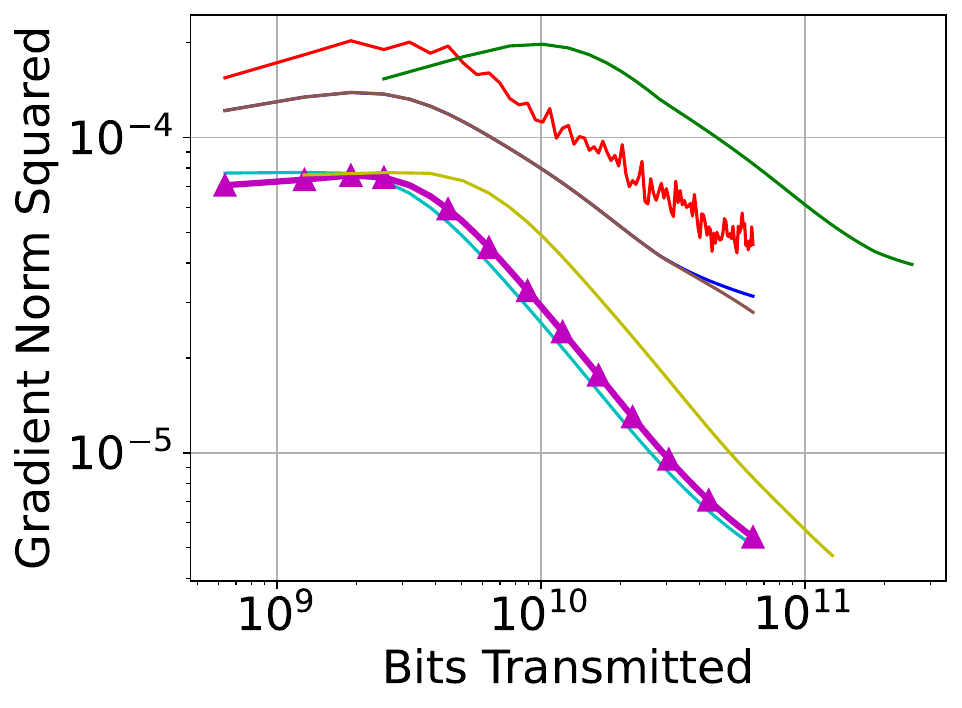} 
    \includegraphics[width=0.235\textwidth]{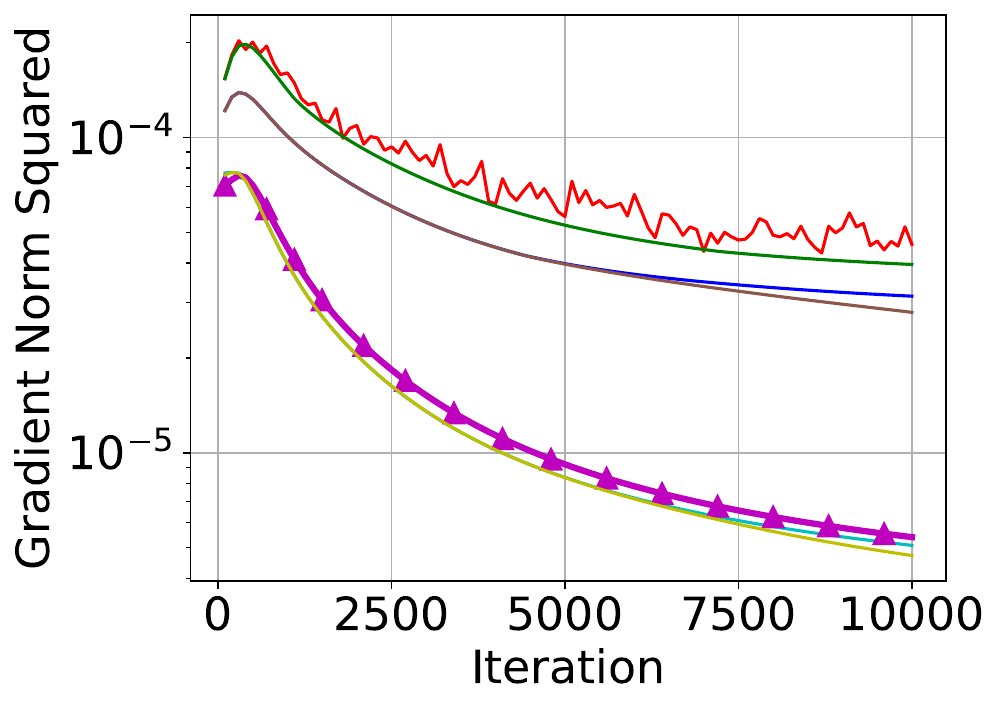} \\
    \includegraphics[width=0.235\textwidth]{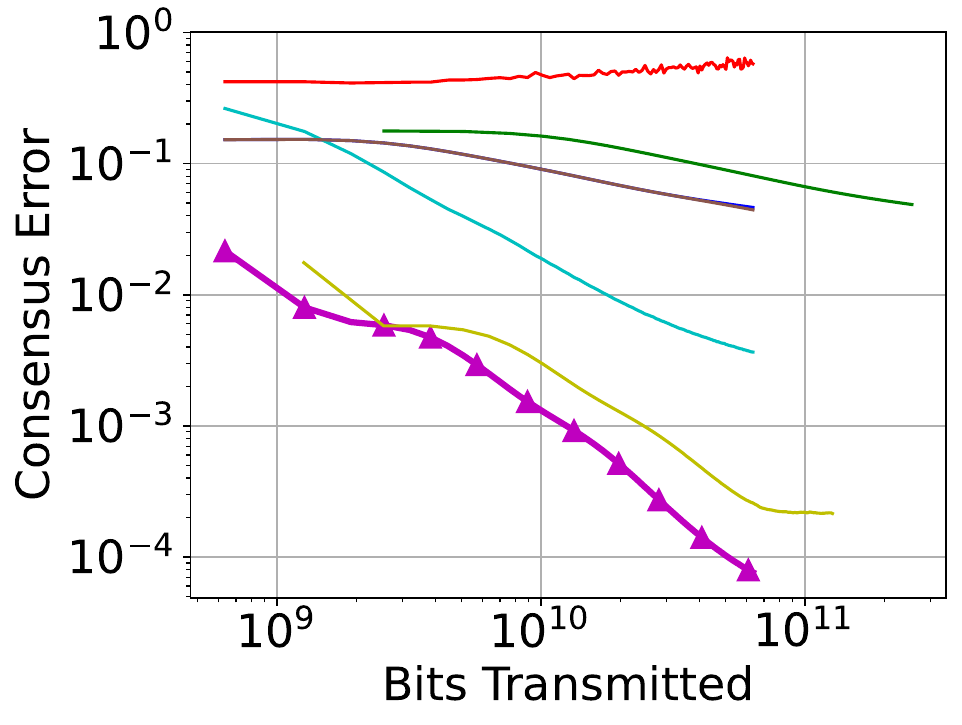}
    \includegraphics[width=0.235\textwidth]{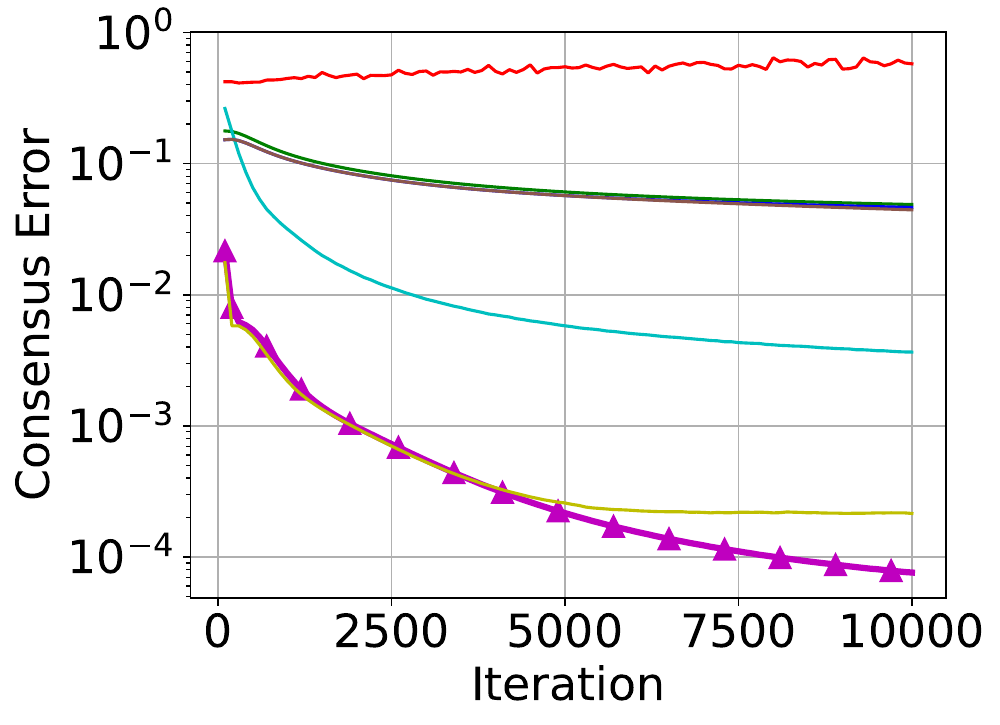} \\
    \includegraphics[width=0.235\textwidth]{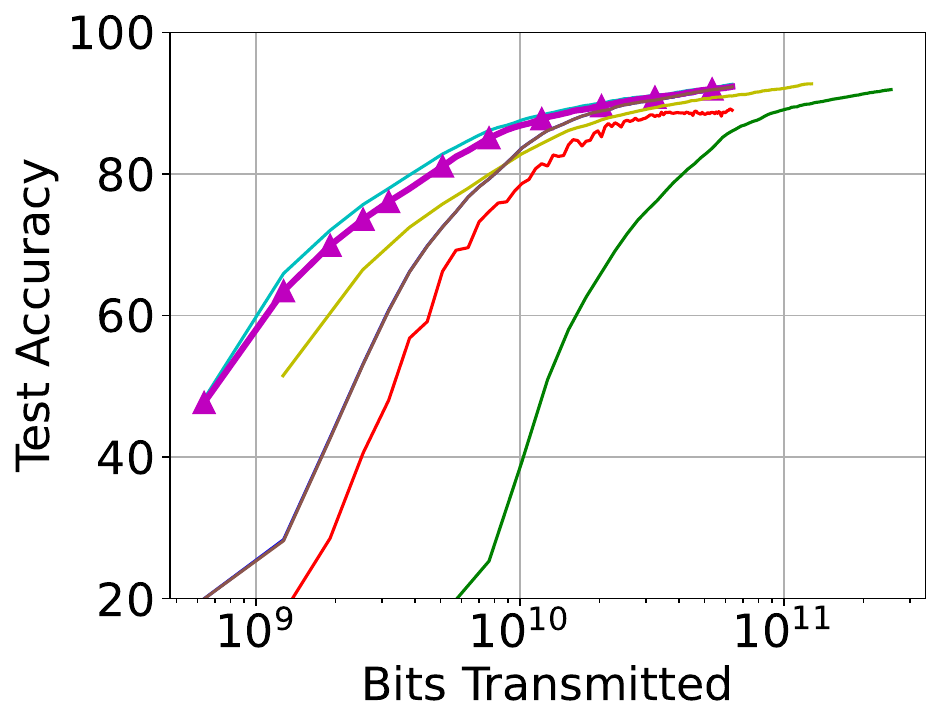} 
    \includegraphics[width=0.235\textwidth]{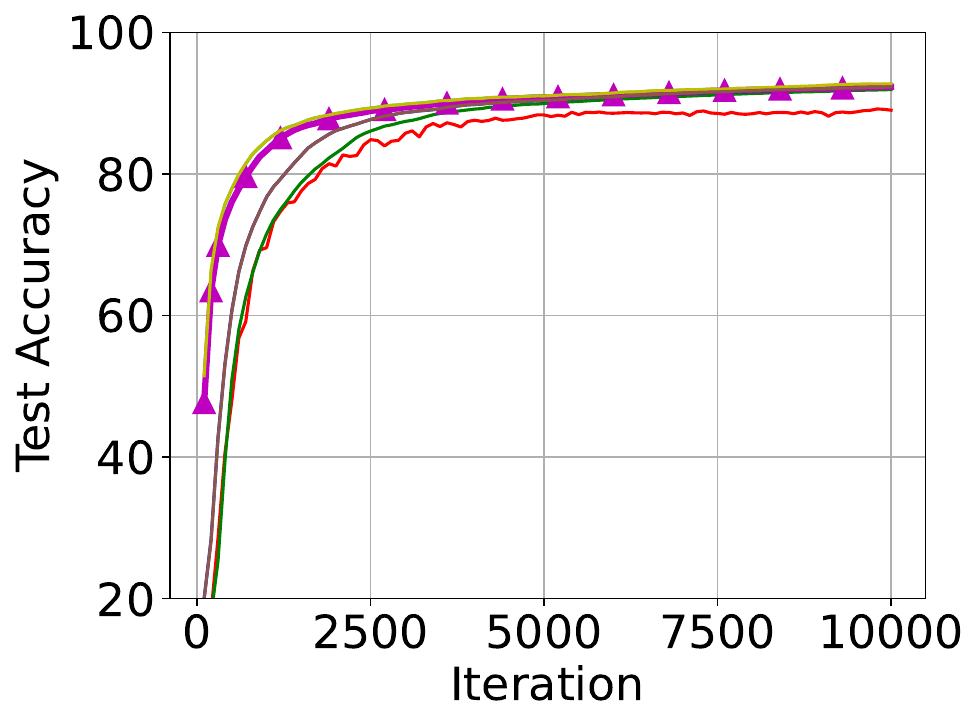} 
    \caption{Training a 2-layer feedforward network using the MNIST data. The bit-rates for communication quantization are displayed in the legend. 
    } 
    \label{fig:mnist}
\end{figure}

Fig.~\ref{fig:mnist} compares the training loss ($\max_{i \in [n]} f(\Prm^t_i)$), gradient norm ($\max_{i\in [n]}\| \nabla f( \Prm^t_i ) \|^2$), consensus error ($\| \Prm^t \|_{ \tilde{\bf K} }^2$), and test accuracy of the worst local model $\Prm^t_i$ against the iteration number and communication cost (in bits transmitted). Notice that for the uncompressed DGD method, we assume that a 32 bits full precision representation for a real number. 
We observe that both the primal-dual algorithms, {\algoname}, CP-SGD, achieve substantially better performance than the DGD, CHOCO-SGD algorithms that are primal-only. It illustrates how data heterogeneity will dampen the convergence of primal-only algorithms. Furthermore, we note that the heuristic modification of DGD with direct quantization results in a non-converging algorithm, indicating the necessity of developing better approaches such as error feedback. When comparing between TiCoPD and CP-SGD, we can see that TiCoPD achieves lower consensus error by two orders of magnitude as the number of iterations grows. Meanwhile, the two algorithms share similar convergence rates on other metrics.


\section{Conclusions}
This paper studies a communication efficient primal-dual algorithm for decentralized optimization with support for compression schemes such as quantized message exchanges. Unlike prior works, our key idea is to develop the algorithm from the augmented Lagrangian framework, and to incorporate classical designs such as majorization-minimization and two-timescale updates. The resultant algorithm is an algorithm that converges at the rate of ${\cal O}(T^{-1})$ for smooth (possibly non-convex) problems. We envisage that the proposed algorithmic framework can be extended for numerous tasks in signal processing and machine learning. 


\newpage

\end{document}